\documentclass[runningheads]{cl2emult}

\usepackage{amsmath,amssymb}

\usepackage{makeidx}  % allows index generation
\usepackage{graphicx} % standard LaTeX graphics tool
\makeindex            % used for the subject index
\usepackage{amsmath}
\begin{document}
\title*{Bounds for Completely Decomposable Jacobians}
\toctitle{Bounds for Completely Decomposable Jacobians}
% allows explicit linebreak for the table of content
%
%
\titlerunning{Bounds for Completely Decomposable Jacobians}
% allows abbreviation of title, if the full title is too long
% to fit in the running head
%
\author{Iwan Duursma\inst{1}
\and Jean-Yves Enjalbert\inst{2}}
\authorrunning{Iwan Duursma et al.}
% if there are more than two authors,
% please abbreviate author list for running head
%
%
\institute{University of Illinois at U-C, Urbana IL 61801, USA
\and Universit\'{e} de Limoges, F-87060 Limoges Cedex, France}

\maketitle              % typesets the title of the contribution

\begin{abstract}
A curve over the field of two elements with completely decomposable Jacobian is shown to have at most 
six rational points and genus at most 26. The bounds are sharp. The previous upper 
bound for the genus was 145.
We also show that a curve over the field of $q$ elements with more than $q^{m/2}+1$ rational points has
at least one Frobenius angle in the open interval $(\pi/m,3\pi/m)$. The proofs make use of the 
explicit formula method.
\end{abstract}

\section{Introduction}
The Jacobian of an (absolutely irreducible, projective, non-singular) algebraic curve is said to be completely 
decomposable if it is isogenous over the base field to a product of elliptic curves. 
Many examples are known of curves with
completely decomposable Jacobian \cite{EkeSe}, both in characteristic zero and in finite characteristic.
For a curve over a finite field $F_q$, the genus of a curve with completely decomposable 
Jacobian is bounded \cite{TsfVl-Asymp}, \cite{Serre-Hecke}. For $q=2$, Serre \cite{Serre-Hecke} gives a
first order estimate $g < 146$. We use the explicit formula method developed in \cite{Serre-Points} to
obtain $g \leq 26$. The upper bound is sharp and is attained by the modular curve $X(11)$ for which Hecke 
showed that the Jacobian
decomposes as $E_1^5 \times E_2^{10} \times E_3^{11}$ \cite{Ligoz}. \\

For an algebraic curve (absolutely irreducible, projective, non-singular) of genus $g$ over a finite field of 
$q$ elements, the Hasse-Weil bound gives that the number
of rational points $N$ does not exceed $q+1+2g\sqrt{q}$.
For the explicit formula method, the number of rational points is expressed in terms of the Frobenius 
eigenvalues as
\[
N = q+1-\sum_{j=1}^{g} (\alpha_j + \bar \alpha_j ).
\]
By Weil's theorem, we may write $\alpha_j = \sqrt{q} e^{i\theta_j}$, for elements $\theta_j$ in $[0,\pi]$ for all $j$.
The $\theta_j$ are called the {\sl Frobenius angles}. Over an extension field of size $q^m$, the number
of rational points $N_m$ is given by
\[
N_m = q^m+1 - \sum_{j=1}^{g} (\alpha_j^m + \bar \alpha_j^m ) = q^m + 1 - r^m \sum_{j=1}^{g} 2 \cos m \theta_j,  
 \]
where $r=\sqrt{q}$.
For curves of large genus, the distribution of the Frobenius angles is restricted by the constraints 
$N_{dm} \geq N_{m}$, for all $m,d$.
This allows one to obtain upper bounds of the form $N \leq ag+b$
for the number of rational points that are better than the Hasse-Weil bound 
when the genus is large \cite{Ihara}, \cite{Serre-Points}.
Asymptotically, the Drinfeld-Vladuts bound gives $\lim \sup_{g \rightarrow \infty} N/g \leq \sqrt{q}-1$ \cite{DriVl},
where the limit is over an infinite family of curves of increasing genus. 
In Section \ref{sec:expl}, we recall the main steps of the explicit formula method. \\

Tsfasman-Vladuts \cite{TsfVl-Asymp} and Serre \cite{Serre-Hecke} study the distribution of the Frobenius angles 
for families of curves of increasing genus. It is
easy to see that any infinite family of curves of increasing genus contains a subfamily for which
$N_m / g$ approaches a limit, for each $m$, when the genus increases. Such subfamilies are called asymptotically 
exact in \cite{TsfVl-Asymp}.
For curves in an asymptotically exact family, the distribution of the Frobenius angles approaches a 
limit distribution that is given by a continuous measure on $[0,\pi]$. In particular, the Frobenius
angles in an asymptotically exact family are dense in $[0,\pi]$. This shows that any family of curves for
which the Frobenius angles are not dense in $[0,\pi]$ is finite. We consider the following problem.

\begin{description}
\item[(Problem 1)] Given a discrete subset $\Theta$ of $[0,\pi]$, 
maximize $N$ and $g$ for a curve over $F_q$ with all Frobenius angles in $\Theta$. 
\end{description}

The elliptic curves over the field of two elements have Frobenius angle 
$\theta$ such that $2 \sqrt{2} \cos \theta \in \{-2, -1, 0, 1, 2\}$. The corresponding Frobenius 
eigenvalues are of degree at most two. As a special case of the previous problem we have

\begin{description}
\item[(Problem 2)] Maximize $N$ and $g$ for a curve over $F_q$ 
with all Frobenius eigenvalues of bounded degree at most $d$. 
\end{description}

The case $d=2$ corresponds to curves with completely decomposable Jacobian. 
In Section \ref{sec:N6} and Section \ref{sec:g26}, respectively, we show that a curve over $F_2$ with completely decomposable
Jacobian has $N \leq 6$ and $g \leq 26$, respectively.
Similarly, the family of curves with no Frobenius angle in a given interval is finite. 
And we can ask for the largest number of rational points or the largest genus for curves in the family. 

\begin{description}
\item[(Problem 3)] Given a (small) subset $I$ of $[0,\pi]$, 
maximize $N$ and $g$ for a curve over $F_q$ with all Frobenius angles outside $I$. 
\end{description}

In Section \ref{sec:asym}, we prove that any curve over ${F}_{q}$ with
$N > q^{m/2}+1$ has a Frobenius angle in the open interval $(\pi/m,3\pi/m)$.
We formulate one other problem along the same lines. It will not be considered in this paper however.

\begin{description}
\item[(Problem 4)] Given $\delta$, maximize $N$ and $g$ for a curve over $F_q$ such that
$[0,\pi] \not \subset \cup_j (\theta_j-\delta,\theta_j+\delta)$. 
\end{description}

\section{The explicit formula method}

We first recall the explicit formula method and its use in obtaining general upper bounds for the
number of rational points on a curve \cite{Serre-Points}, \cite{Hansor}. Then we present three
variations of the method that yield better bounds for curves whose Frobenius angles are restricted
to a subset $\Theta$ of $[0,\pi]$. In particular, curves that exceed one of the latter bounds, 
necessarily have at least one Frobenius angle outside $\Theta$.  

\subsection{General upper bounds for the number of rational points} \label{sec:expl}

For an algebraic curve $X$ of genus $g$ over the finite field $F_q$ of $q$ elements, 
let the Frobenius angles be $\theta_1, \theta_2, \ldots, \theta_g$. So that 
the number of rational points $N_n$ over $F_{q^n}$ satisfies
\[
N_n  = q^n + 1 - q^{n/2} \sum_{j=1}^g 2 \cos n \theta_j.
\]
With $r = \sqrt{q}$,
we rewrite the equation as
\begin{equation} \label{eq:Nn}
N_1 r^{-n} + (N_n-N_1) r^{-n} = ~~r^n + r^{-n} - \sum_{j=1}^g 2 \cos n \theta_j. 
\end{equation}
Let $f$ be an auxiliary cosine polynomial with real coefficients $u_n$,
\begin{equation} \label{eq:f}
f(\theta ) = u_0 + \sum_{n \geq 1} u_n \cos n \theta.
\end{equation}
Define 
\begin{equation} \label{eq:psi}
\psi(x) = \sum_{n \geq 1} u_n x^n.
\end{equation}
The equations (\ref{eq:Nn}) scaled by $u_n$, for $n=1,2,\ldots$, add up to
\begin{multline} \label{eq:exp}
N_1 \psi (r^{-1}) + \sum_{n \geq 2} u_n (N_n-N_1) r^{-n} = \\ 
= 2 u_0 g + \psi(r) + \psi(r^{-1})  - 2 \sum_{j=1}^{g} f(\theta_j).
\end{multline}
The equation (\ref{eq:exp}) leads to upper bounds for the number of points.
As in \cite{Serre-Points}, choose $\{ u_n \}$ such that $u_0 = 1$, and
\begin{description}
\item[(a)] ${u}_{n}\geq 0, \forall n \geq 1$ 
\item[(b)] $f(\theta)\geq 0, \text{for all } \theta \in [0,\pi].$
\end{description}
Then Equation (\ref{eq:exp}) yields
\[
N \psi(r^{-1}) \leq 2g + \psi(r^{-1}) + \psi(r).
\]
As an example, the choice
\begin{align*}
f (\theta ) &= \cos^2 \theta (1 - \cos \theta / \cos (\frac{5\pi}{6}))^2 \\
  &= 1 + \sqrt{3} \cos \theta +  \frac{7}{6} \cos 2 \theta 
     + \frac{\sqrt{3}}{3} \cos 3 \theta +  \frac{1}{6} \cos 4 \theta
\end{align*}
gives, for $q=3$, the upper bound
\[
N \leq \frac{54}{41} (g-15) + 28 < 1.317 \; g + 8.244.
\]
This is better than the Hasse-Weil bound $N \leq 2 \sqrt{3} \, g + 4$ for all $g \geq 2$. 
A curve attains the upper
bound above only if $N_1 = N_2 = N_3 = N_4$ and if all its Frobenius angles are among $\{ \pi/2, 5\pi/6 \}$. 
The unique such curve is the Deligne-Lusztig curve associated to ${}^2 G_2(3)$ \cite{HanPe}. 
The curve is of genus $g=15$ and has $N=28$. 
Its zeta function $Z(T)=P(T)/(1-T)(1-3T)$ has numerator $P(T) = (1+3T^2)^7 (1+3T+3T^2)^8$.
% In particular, the curve has completely decomposable Jacobian.

\subsection{Restricted upper bounds for the number of rational points ($u_0=1$)} \label{sec:angfrob}

The upper bound in the previous subsection generalizes as follows. Choose $\{ u_n \}$
in Equation (\ref{eq:exp}) such that
\begin{description}
\item[(a)] $u_0 = 1$ and ${u}_{n}\geq 0, \forall n \geq 2$. 
\item[(b)] $f(\theta)\geq 0, \text{for all } \theta \in \Theta \subset [0,\pi].$
\end{description}
Then, for a curve that has all its Frobenius angles contained in $\Theta$,
% for (\ref{eq:exp}) we know that every curves with
\[
N \psi(r^{-1}) \leq 2g + \psi(r^{-1}) + \psi(r).
\]
The converse yields that a curve with
\[
N \psi(r^{-1}) > 2g + \psi(r^{-1}) + \psi(r).
\]
has a Frobenius angle outside $\Theta$. For $0 < \alpha < \beta < \pi$, let 
\begin{align*}
f_2(\theta) &= \left( \cos \theta - \cos \alpha \right) \left( \cos \theta - \cos \beta \right), \\
          &= \frac{1}{2}+\cos \alpha \cos \beta - \left( \cos \alpha + \cos \beta \right) \cos \theta
               + \frac{1}{2}\cos 2\theta.
\end{align*}
Then $f_2(\theta)$ is non-negative on $\Theta = [0,\pi] \backslash (\alpha,\beta)$.
For $q=2$, and for $\alpha=\pi/3$ and $\beta=3\pi/4$, we obtain
\[
N > \frac{8-2 \sqrt{2}}{7} (g-1) + 5  ~~\Rightarrow~~ \exists \theta_j \in (\frac{\pi}{3},\frac{3\pi}{4}).
\] 
The inequality on the left applies in the range $2 \leq g \leq 38$. In that range the inequality holds for a curve 
that meets the Oesterl\'e upper bound for the number of points.  
For another example, let
\begin{align*}
f(\theta) &= (1+\sqrt{2} \cos \theta)(1-2\sqrt{2} \cos \theta)^2, \\
          &= 1 + 3 \sqrt{2} \cos \theta + 2 \sqrt{2} \cos 3 \theta.
\end{align*}
We obtain, for a curve over $F_2$,
\[
N > \frac{1}{2} (g-1) + 5 ~~\Rightarrow~~ \exists \theta_j \in (\frac{3\pi}{4},\pi].
\] 

\subsection{Uniform upper bounds for the number of rational points ($u_0=0$)} \label{sec:N6}
 
By choosing $u_0 = 0$, we obtain upper bounds for the number of rational points that are
independent of the genus $g$. Choose $\{ u_n \}$ in Equation (\ref{eq:exp}) such that
\begin{description}
\item[(a)] $u_0 = 0$ and ${u}_{n}\geq 0, \forall n \geq 2$. 
\item[(b)] $f(\theta)\geq 0, \text{for all } \theta \in \Theta \subset [0,\pi].$
\end{description}
Then the number $N$ of rational points on a curve with all Frobenius angles contained in $\Theta$ satisfies 
\[
N \psi(r^{-1}) \leq \psi(r^{-1}) + \psi(r).
\]
If, moreover, the coefficients $u_n$ have the following symmetry property, for some positive integer 
$m$ with $m>\deg(\psi)$, 
\begin{description}
\item[(c)] $u_{n} = u_{m-n}$, for $n=0,1,\ldots,m,$
\end{description}
then the upper bound becomes
\[
N \leq 1+ \frac{\sum_{n = 0}^m u_n r^n}{\sum_{n = 0}^m u_{m-n} r^{n-m}} = r^m + 1.
\]
The function
\begin{align*}
f(\theta ) &= \frac{\sqrt{2}}{5} \cos \theta (1 - 2 \cos^2 \theta) (1 - 8 \cos^2 \theta) \\
  &= {\frac{7}{10}}\,\sqrt {2}\cos \theta + {\frac {1}{2}}\,\sqrt {2}\cos 3\,\theta 
       + {\frac {1}{5}}\,\sqrt {2}\cos 5\,\theta 
\end{align*}
cancels at the Frobenius angles of the five different elliptic curves over $F_2$.
It leads to the bound $N \leq 6$ for any curve $X$ over $F_2$ with completely decomposable Jacobian.
The bound is tight only when $N_1 = N_3 = N_5$.
The smallest feasible zeta function is of genus $3$ with uniquely determined zeta polynomial
$P(T)=(1+2T+2T^2)^2{\linebreak}(1-T+2T^2).$ It is realized by the curve 
\[
y^2+y = \frac{x^2+x}{(x^2+x+1)^3}.
\]
We give two examples that use Condition (c). The choice $f(\theta) = \cos \theta$ yields that 
a curve with $N > r^2+1$ has a Frobenius angle in $(\pi/2,3\pi/2)$ (indeed
the Frobenius trace can only be negative if at least one Frobenius angle has $\cos \theta < 0$).
The choice $f(\theta) = \cos \theta + \cos 2\theta$ yields that 
a curve with $N > r^3+1$ has a Frobenius angle in $(\pi/3,\pi)$. In both cases, the bound on $N$ 
is sharp. The projective line with $N=r^2+1$ has no Frobenius angle in $(\pi/2,3\pi/2)$, and the
Hermitian curve (see \cite{StiRu}) over $F_{r^2}$ with $N=r^3+1$ has no Frobenius angle in $(\pi/3,\pi)$. 
The latter
example confirms that the Hasse-Weil bound is not sharp for curves with $N > r^3+1$.
In Section \ref{sec:asym}, we show more generally that a curve with $N > r^m+1$ has a Frobenius
angle in $(\pi/m,3\pi/m)$. 
 
\subsection{Uniform upper bounds for the genus ($u_0=-1$)} \label{sec:g26}

By choosing $u_0 = -1$, we obtain upper bounds for the genus $g$. 
Choose $\{ u_n \}$ in Equation (\ref{eq:exp}) such that
\begin{description}
\item[(a)] $u_0 = -1$ and ${u}_{n}\geq 0, \forall n \geq 2$. 
\item[(b)] $f(\theta)\geq 0, \text{for all } \theta \in \Theta \subset [0,\pi].$
\end{description}
Then the genus of a curve with all Frobenius angles contained in $\Theta$ satisfies
$$N \psi(r^{-1}) + 2g \leq \psi (r) + \psi (r^{-1}).$$
If, moreover, the coefficients $u_n$ satisfy
\begin{description}
\item[(d)] $\psi(r^{-1})=0,$
\end{description}
then the upper bound becomes
$$2g \leq \psi (r). $$
The function 
\[
f(\theta) = {-1-\frac {4}{3}}\,\cos \theta +{\frac {7}{9}}\,\cos 2\,\theta \,
+{\frac {26}{9}}\,\cos 3\,\theta +{\frac {16}{9}}\,\cos 4\,\theta 
\]
is of minimal degree such that it cancels at the Frobenius angles of the three elliptic curves over 
$F_4$ that are defined over $F_2$ and such that Condition (d) holds.
It leads to the bound $2g \leq 52$ for any curve $X$ over $F_2$ with completely decomposable Jacobian. 
A previous estimate showed that $g \leq 145$ \cite{Serre-Hecke}.
The bound is tight only when $N_1 = N_2 = N_3 = N_4$ for the base field $F_4$. 
It is attained by the modular curve $X(11)$, which has $g=26$, $N=55$ over $F_4$, 
and zeta polynomial $P(T)=(1+4T+4T^2)^5(1+3T+4T^2)^{10}(1+4T^2)^{11}$. \\

\section{An asymptotic example} \label{sec:asym}

Let $m \geq 4$ and let $\alpha = \pi/m$. Conditions (a)-(c) in Section \ref{sec:N6} hold with 
$\Theta = [0,\pi] \backslash (\pi/m,3\pi/m)$
for coefficients $\{u_n\}$ that are defined by
\begin{equation} \label{eq:ftheta}
f(\theta) = \frac{1+\cos m \theta }{4(\cos \theta - \cos \alpha)(\cos \theta - \cos 3\alpha)} =
\sum_{n=2}^{m-2} u_n \cos n \theta.  
\end{equation}
So that
\begin{equation}  \label{eq:uj}
u_n =  \frac{\sin (n-1) \alpha \, \sin n \alpha \, \sin (n+1) \alpha}
       {\sin \alpha  \sin 2 \alpha  \sin 3 \alpha}, \quad n=0,1,\ldots,m.  
\end{equation}
Thus, a curve with number of rational points $N>{r}^{m}+1$, for $m \geq 4$, has at least one Frobenius angle in the open interval
$(\pi/m,3\pi/m)$. For $f(\theta)$ we may write
\[  
f(\theta) = 2^{m-3} \prod_{k=2}^{m-1}\big( \cos \theta - \cos(2k+1)\alpha \big),
\]
which justifies writing the right hand side of (\ref{eq:ftheta}) as a cosine polynomial. 
To see that the coefficients of the cosine polynomial are those given by (\ref{eq:uj}), 
we use a generating function for gaussian polynomials \cite{AndrewsReprint}
\[
\frac{1}{(1-T)(1-yT)(1-y^2T)(1-y^3T)} = \sum_{i\geq0} {\left[ \begin{array}{c}i+3\\ 3\end{array} \right]} T^i,
\]
where 
\[
{\left[ \begin{array}{c}i+3\\ 3\end{array} \right]} = \frac{(y^{i+3}-1)(y^{i+2}-1)(y^{i+1}-1)}{(y^3-1)(y^2-1)(y-1)}.
\]
For $y$ with $y^m=1$, the right hand side is periodic and, for $n=i+2$,
\[
\frac{T^2(1-T^m)}{(1-T)(1-yT)(1-y^2T)(1-y^3T)} = \sum_{n=2}^{m-2} \frac{(y^{n+1}-1)(y^{n}-1)(y^{n-1}-1)}{(y^3-1)(y^2-1)(y-1)} T^n.
\]
Let $x=e^{i\alpha}$, so that $x^m=-1$. With $y=x^2$ and $t=x^3T$, we obtain
\[
\frac{(1+t^m)}{(t+t^{-1}-2\cos \alpha)(t+t^{-1}-2\cos 3\alpha)} = \sum_{n=2}^{m-2} 
\frac{\sin (n-1) \alpha \, \sin n \alpha \, \sin (n+1) \alpha}
       {\sin \alpha  \sin 2 \alpha  \sin 3 \alpha} t^n.
\]
Now sum the two equations with $t=e^{i\theta}$ and $t=e^{-i\theta}$, respectively, and divide
by 2. \\

The cases $m=2$ and $m=3$ were considered in Section \ref{sec:N6}, so that the claim extends to all $m \geq 2$.
For $m=4$ and $m=6$ the bounds are sharp, as can be seen by considering curves of
Suzuky type or Ree type, respectively. The Suzuki curve over $F_8$ has $N=65$ but has no Frobenius angle in $(\pi/4,3\pi/4)$.
The Ree curve over $F_3$ has $N=28$ but has no Frobenius angle in $(\pi/6,\pi/2)$.

\section{Conclusion}

Results by Tsfasman-Vladuts and Serre led us to consider Problems (1)-(4) in the Introduction. 
For Problems (1)-(3) we have given methods that yield partial results. One result is a sharp upper bound
for the number of points ($N \leq 6$) or the genus ($g \leq 26$) for a curve over $F_2$ with completely
decomposable Jacobian. We also showed that a curve over $F_q$ with $N > q^{m/2}+1$ has at least one Frobenius angle in the
interval $(\pi/m,3\pi/m)$. No results were obtained towards Problem (4).  

\nocite{EkeSe}
\nocite{HanPe}
\nocite{Ligoz}
\nocite{Serre-Points}
\nocite{Serre-Hecke}
\nocite{StiRu}
\nocite{Tsfas-Asymp}
\nocite{TsfVl-Asymp}
\nocite{Hansor}
%\nocite{AndrewsOriginal}
\nocite{AndrewsReprint}

\def\cprime{$'$}

%\bibliographystyle{alpha}

%\bibliography{Fq6}

\newpage

\section*{Addendum : The modular curve $X(11)$ modulo $2$}

\newcommand{\z}{\zeta}
\newcommand{\bt}{{\underline \beta}}
\newcommand{\ch}{{\rm char}}

The main text uses the following claim. \\

Claim: The modular curve $X(11)$ of genus $26$ has
a model over ${\mathbb F}_2$ with $55$ rational points over ${\mathbb F}_4$. \\

The claim follows from results by Klein together with the following fact about the
group $PSL(2,11)$. \\

Fact 1: The subgroups of order $12$ in $PSL(2,11)$ divide into two conjugacy classes, each
consisting of $55$ copies of $A_4$. In each case the group $A_4$ is normally closed
in $PSL(2,11)$. \\

The $55$ rational points over ${\mathbb F}_4$ on $X(11)$ are the points above $j=0$ in
the galois cover $X(11) \longrightarrow X(1)$ of degree $660$ with group $PSL(2,11)$.
Let the cover be
defined over a finite field $k$ of characteristic $2$ such that all automorphisms
have their coefficients in $k$. Fact 1 shows that the stabilizer of a point above $j=0$
is a subgroup $A_4$ and moreover, since $A_4$ is normally closed, that the point is
the unique fixed point of the stabilizer. In particular, the points above $j=0$ are
defined over $k$. \\

With the previous argument it suffices to find a model for $X(11)$ defined over 
${\mathbb F}_2$ such that the automorphisms are defined over ${\mathbb F}_4$.
As starting point, we use Klein's model for $X(11)$ in ${\mathbb P}^4({\mathbb C})$. \\
% defined over ${\mathbb Q}$. \\
%with the following properties. \\

Fact 2 [Klein, Ges. Math. Abh. III, p.146]: Klein's model for $X(11)$ is defined over 
${\mathbb Q}$ and has all its $660$ automorphisms defined over ${\mathbb Q}(\rho)$, 
where $\rho$ is a primitive $11$-th root of unity. The group contains
the cyclic automorphism $C(v:w:x:y:z)=(z:v:w:x:y)$. Let $F$ be the automorphism of 
${\mathbb Q}(\rho)$ with $F(\rho) = \rho^4$. The generators $S$ (of order $11)$ and $T$ 
(of order $2$) presented by Klein satisfy $F(S) = C^{-1} \circ S \circ C$ and
$F(T) = C^{-1} \circ T \circ C$. \\

%While $S$ and $T$ generate the group, Klein includes $C$ as a third generator so that the
%elements of the automorphism group can be described as $C^\alpha S^\beta$, $C^\alpha S^\beta T S^\gamma$, 
%for $\alpha = 0,1,2,3,4$ and $\beta, \gamma = 0,1,2, \ldots, 10$. \\
A change of variables $(v:w:x:y:z)=Q(v':w':x':y':z')$ leads to a 
new model in the variables $(v':w':x':y':z')$. The new model is defined over ${\mathbb Q}$
if the set $\{ Q, C \circ Q, C^2 \circ Q, C^3 \circ Q, C^4 \circ Q \}$ is defined over 
${\mathbb Q}$. This is achieved by defining $Q$ over ${\mathbb Q}(\rho+\rho^{-1})$ such 
that $F(Q) = C^i \circ Q$, for some $i \in \{0,1,2,3,4\}$. After the change of variables, 
the automorphims are generated by $Q^{-1} \circ S \circ Q$ and $Q^{-1} \circ T \circ Q$.  
For a choice of $Q$ with $F(Q) = C^{-1} \circ Q$, the automorphisms are defined over 
the fixed field of $F$: $F(Q^{-1} \circ A \circ Q) = Q^{-1} \circ A \circ Q,$ for $A = S, T.$ \\
%\begin{align*}
%F(Q^{-1} \circ S \circ Q) &= F(Q^{-1}) \circ F(S) \circ F(Q) \\
% &= (Q^{-1} \circ C) \circ ( C^{-1} \circ S \circ C) \circ (C^{-1} \circ Q) \\
%  &= Q^{-1} \circ S \circ Q.
%\end{align*} 
%If $f(Av)=f(v)$ and $g(v')=f(Qv')=0$ then $g(Bv')=g(v')$ if $f(QBv') = f(Qv')$ if $QB=AQ$.

Klein's model has good reduction modulo $2$ and the twisted model is defined over ${\mathbb F}_2$ 
such that all automorphisms are defined over ${\mathbb F}_4$.

\end{document}